\documentclass[a4paper,12pt]{amsart}
\usepackage{amssymb, amsmath}
\usepackage{amscd}
\numberwithin{equation}{section}
\usepackage{epsfig}
\usepackage{amsmath}
\usepackage{amsfonts,amssymb,amsopn}
\usepackage[all]{xy}
\SilentMatrices

\newcommand{\pp}{\mathbb{P}}
\newcommand{\cc}{\mathbb{C}}
\newcommand{\zz}{\mathbb{Z}}
\newcommand{\A}{\mathbb{A}}

\newcommand{\Homom}{\mathrm{Hom}}
\newcommand{\Endom}{\mathrm{End}}
\newcommand{\Ker}{\mathrm{Ker}}

\newcommand{\Image}{\mathrm{Im}}

\newcommand{\rank}{\mathrm{rk}}
\newcommand{\Sym}{\mathrm{Sym}}

\newcommand{\Ox}{\mathcal{O}_{X}}
\newcommand{\Kx}{K_{X}}
\newcommand{\sRat}{\underline{\mathrm{Rat}}}
\newcommand{\Rat}{\mathrm{Rat}}
\newcommand{\sPrin}{\mathrm{\underline{Prin}}}
\newcommand{\Prin}{\mathrm{Prin}}

\newcommand{\Sec}{\mathrm{Sec}}

\newcommand{\MO}{\mathcal MO}
\newcommand{\MS}{\mathcal MS}

\newcommand{\M}{\mathcal M}

\newcommand{\cV}{\mathcal{V}}
\newcommand{\cW}{\mathcal{W}}
\newcommand{\Iden}{\mathrm{Id}}
\newcommand{\GL}{\mathrm{GL}}
\newcommand{\Sp}{\mathrm{Sp}}
\newcommand{\SO}{\mathrm{SO}}
\newcommand{\SL}{\mathrm{SL}}
\newcommand{\Gr}{\mathrm{Gr}}
\newcommand{\sy}{\mathbb{S}}
\newcommand{\tU}{\tilde{U}}
\newcommand{\cE}{\mathcal{E}}
\newcommand{\e}{\varepsilon}

\newtheorem{theorem}{{\textbf Theorem}}[section]
\newtheorem{proposition}[theorem]{{\textbf Proposition}}
\newtheorem{corollary}[theorem]{{\textbf Corollary}}
\newtheorem{lemma}[theorem]{{\textbf Lemma}}
\newtheorem{criterion}[theorem]{{\textbf Criterion}}
\newtheorem{defn}[theorem]{{\textbf Definition}}

\newtheorem{remit}[theorem]{{\textbf Remark}}
\newenvironment{remark}{\begin{remit}\rm}{\end{remit}}
\newenvironment{definition}{\begin{defn}\rm}{\end{defn}}

\title[symplectic and orthogonal bundles]{A stratification on the moduli spaces of \\  symplectic  and orthogonal bundles over a curve}

\author{Insong Choe and George H.\ Hitching}

\begin{document}

\begin{abstract}
A symplectic or orthogonal bundle $V$ of rank $2n$ over a curve has an invariant $t(V)$ which measures the maximal degree of its isotropic subbundles of rank $n$. This invariant $t$ defines stratifications on moduli spaces of symplectic and orthogonal bundles. We study this stratification by relating it to another one given by secant varieties in certain extension spaces.

We give a sharp upper bound on $t(V)$, which generalizes the classical Nagata bound for ruled surfaces and the Hirschowitz bound for vector bundles, and study the structure of the stratifications on the moduli spaces. In particular, we compute the dimension of each stratum. We give a geometric interpretation of the number of maximal Lagrangian subbundles of a general symplectic bundle, when this is finite. We also observe some interesting features of orthogonal bundles which do not arise for symplectic bundles, essentially due to the richer topological structure of the moduli space in the orthogonal case.
\end{abstract}

\maketitle

\section{Introduction}

Let $X$ be a smooth algebraic curve over $\cc$. A vector bundle $V$ of rank two over $X$ determines a ruled surface $\pp V$. Such surfaces have been studied since the 19th century. A line subbundle $L$ of $V$ gives a section $\sigma_L$ of $\pp V$, called a \textsl{directrix curve}. The self-intersection number of $\sigma_L$ is given by
\[
\sigma_{L} \cdot \sigma_{L}  \ = \ \deg (V/L) - \deg L \ = \ \deg V - 2 \deg L.
\]
The \textsl{Segre invariant} of $\pp V$ is defined as the minimal value of $\sigma_{L} \cdot \sigma_L$ over all $L \subset V$. Via the above formula, this invariant also provides a measure of the difference of the slopes of $V$ and $L$.

The Segre invariant yields a natural stratification on the moduli space of vector bundles of rank two over $X$, which was studied by Lange and Narasimhan \cite{LN}. A generalization of this stratification to the moduli of vector bundles over $X$ of arbitrary rank was considered by Lange \cite{Lan} and the details were settled in Brambila-Paz--Lange \cite{BPL}, and Russo--Teixidor i Bigas \cite{RT}. The aim of the present article is to establish parallel results for symplectic and orthogonal bundles over $X$.

Let us  define a few notions and fix notations. A vector bundle $V$ over $X$ of rank $2n$ is called \textsl{symplectic} (resp.,\ \textsl{orthogonal}) if there is a nondegenerate alternating (resp.,\ symmetric) bilinear form $\omega \colon V \otimes V \to \Ox$. A subbundle $E$ of $V$ is called \textsl{isotropic} if $\omega|_{E \otimes E} = 0$. By linear algebra, an isotropic subbundle of $V$ has rank $\le n$. When $V$ is symplectic, a rank $n$ isotropic subbundle of $V$ is often called a \textsl{Lagrangian subbundle}. We say that a symplectic or orthogonal bundle $V$ is \textsl{stable} (resp., \textsl{semistable}) if for every isotropic subbundle $E$ of $V$, we have
\[
\mu(V) = \frac{ \deg V}{\rank V} \ >  \ \frac{ \deg  E}{\rank E} = \mu (E) \ \  \ (\text{resp., } \mu(V) \ge \mu(E)).
\]
Note that this is a priori weaker than the stability condition for $V$ as a vector bundle; compare with Ramanathan \cite{Rth}. However, Ramanan \cite{Ram} proved that semistability as an orthogonal bundle is equivalent to semistability of the underlying vector bundle, and moreover that a general stable orthogonal bundle is a stable vector bundle. The same argument (worked through in \cite{Hit0}) shows that the analogous statement is true for symplectic bundles.

We denote by $SU_X(2n,\Ox)$ the moduli space of semistable vector bundles of rank $2n$ and trivial determinant, and write $\MS_X(2n)$ (resp., $\MO_X(2n)$) for the sublocus in $SU_X(2n,\Ox)$ of bundles admitting a symplectic (resp., orthogonal) structure.

In the symplectic case it has been proven by Serman \cite{Ser} that the forgetful map 
\[\M_X(\Sp_{2n}\cc) \to SU_X(2n, \Ox)\]
associated to the extension of the structure group $\Sp_{2n}\cc \subset \SL_{2n} \cc$, is an embedding, where $\M_X(\Sp_{2n}\cc)$ is the moduli space of semistable principal $\Sp_{2n}\cc$-bundles over $X$. So $\MS_X(2n)$ coincides with the embedded image of $\M_X(\Sp_{2n}\cc)$.

The orthogonal case is more delicate. By \cite{Ser}, the forgetful map 
\[\M_X(\SO_{2n}\cc) \to SU_X(2n, \Ox)\]
is generically two-to-one, amounting to forgetting the data of an orientation on a principal $\SO_{2n}\cc$-bundle. On the other hand, the map
\[\M_X(\mathrm{O}_{2n}\cc) \to \M_X(\GL_{2n}\cc)\]
is an embedding. The moduli space $\M_X(\mathrm{O}_{2n}\cc)$ of semistable principal $\mathrm{O}_{2n}\cc$-bundles over $X$ has several components, which are indexed by the first and second Stiefel--Whitney classes $(w_{1}, w_{2}) \in H^{1}(X, \zz_{2}) \times H^{2}(X, \zz_{2})$. The class $w_1$ corresponds to the determinant, and there are two components of $\M_X(\mathrm{O}_{2n}\cc)$ with $w_1$ trivial. We write $\MO_{X}(2n)^\pm$ for the embedded images of these components in $SU_X(2n, \Ox)$. 

Let $V$ be a symplectic or orthogonal bundle of rank $2n$. Generalizing the Segre invariant for bundles of rank two, we define
\[
t(V) := -2\max \{  \deg E  : E \hbox{ a rank $n$ isotropic subbundle of } V \}.
\]
In particular, if $V$ is stable (resp.,\ semistable), then $t(V) > 0$ (resp., $t(V) \ge 0$). (Note that when $V$ is a symplectic bundle of rank $>2$, this differs from the invariant $s_{\mathrm{Lag}}$ defined in \cite{CH2} by a constant: $t(V) = \frac{2}{n+1} s_{\mathrm{Lag}}(V) $.) For a symplectic or orthogonal bundle $V$, we denote by $M(V)$ the space of rank $n$ isotropic subbundles $E \subset V$ such that $t(V) = -2 \deg E$. It can be viewed as a closed subscheme of a Quot scheme, so it has a natural structure of projective variety.

The invariant $t(V)$ induces stratifications on moduli spaces of semistable symplectic and orthogonal bundles over $X$.   
For each positive even integer $t$, we define
\[
\MS_X(2n;t)  :=  \{ V \in \MS_X(2n) : t(V) = t \}
\]
and
\[
\MO_X(2n;t)  :=  \{ V \in \MO_X(2n) : t(V) = t \} .
\]
By semi-continuity of the invariant $t(V)$, these subloci are constructible sets.
For the symplectic case, we show:

\begin{theorem} \label{mainsymp}
 \ Let $X$ be a smooth algebraic curve of genus $g \ge 2$.
\begin{enumerate}
\renewcommand{\labelenumi}{(\arabic{enumi})}
\item For any symplectic bundle $V$ of rank $2n$, we have $t(V) \le n(g-1) +1$. This is sharp in the sense that for a general $V \in \MS_X(2n)$, 
\[
n(g-1) \le t(V)  \le n(g-1 ) +1.
\]
\item For each even integer $t$ with $ 2 \le t \le n (g-1)$, the stratum $\MS_X(2n;t)$ is nonempty and irreducible, and
\[
\dim \MS_X(2n;t) \ = \  \frac{1}{2} \left( n(3n+1)(g-1) + (n+1)t \right).
\]
Furthermore, if $t < n(g-1)$, then $\MS_X(2n;t)$ is contained in the closure of $\MS_X(2n;t+2)$ in $\MS_X(2n)$. \\
\item For each positive even integer  $t < n(g-1)$, the space $M(V)$ is a single point for a general $V \in \MS_X(2n;t)$. On the other hand, for a general $V \in \MS_X(2n)$, 
 \[
\dim M(V) \ = \ \begin{cases} 
0&  \text{when} \  n(g-1) \text{ is even, }  \\
\frac{n+1}{2}&  \text{when} \  n(g-1) \text{ is odd}.
\end{cases}
\]
\end{enumerate}
\end{theorem}

Let us give some historical remarks. For vector bundles, the Segre stratifications are well understood; see Lange--Narasimhan \cite{LN}, Hirschowitz \cite{Hir}, Brambila-Paz--Lange \cite{BPL} and Russo--Teixidor i Bigas \cite{RT}. Holla and Narasimhan \cite{HN} defined a generalized Segre invariant for principal $G$-bundles for an arbitrary reductive group $G$, and obtained a bound on the invariant in general. This bound was sharpened for symplectic bundles ($G = \Sp_{2n}\cc$) by the present authors \cite{CH2}, where Theorem \ref{mainsymp} (1) was proven using the Terracini lemma in projective geometry. Also   Theorem \ref{mainsymp} was  proven for symplectic bundles of rank four in \cite{CH2}. The special case when the curve has genus two has earlier been studied in detail  by the second named author \cite{Hit2}. In  this paper, we provide another simpler proof of (1), and prove (2) and (3) for arbitrary rank $n$ and genus $g$. 

In the orthogonal case, the moduli space $\MO_X(2n)$ has two connected components, as was discussed before. We denote by $\MO_{X}(2n)^+$ (resp. $\MO_X(2n)^-$) the component consisting of bundles of trivial (resp., nontrivial) second Stiefel--Whitney class.  We first prove the following.
 
\begin{theorem} \begin{enumerate}
\renewcommand{\labelenumi}{(\arabic{enumi})}
\item Let $E_1$ and $E_2$ be  isotropic  rank $n$ subbundles of an orthogonal bundle  of rank $2n$. Then $\deg E_1$ and $\deg E_2$ have the same parity. \\
\item A semistable orthogonal bundle $V$ belongs to $\MO_{X}(2n)^+$ (resp., $\MO_{X}(2n)^-$) if and only if its isotropic rank $n$ subbundles have even degree (resp., odd degree).
\end{enumerate} \end{theorem}

We then show the following on the stratification on each component.

\begin{theorem} \label{mainorth} \ Let $X$ be a smooth algebraic curve of genus $g \ge 2$. 
\begin{enumerate}
\renewcommand{\labelenumi}{(\arabic{enumi})}
\item For any orthogonal bundle $V$ of rank $2n$, we have $t(V) \le n(g-1) + 3$. This is sharp in the sense that two even numbers $t$ with $n(g-1) \le t \le n(g-1)+3$ correspond to the values of $t(V)$ for general $V$ in $\MO_X(2n)^+$ and $\MO_X(2n)^-$.
\item For each even integer $t$ with $2 \le t  \le n (g-1) $, the stratum $\MO_X(2n;t)$ is nonempty and irreducible, and
\[
\dim \MO_X(2n;t) = \frac{1}{2} \left( n(3n-1)(g-1) + (n-1)t \right).
\]
Furthermore, if $t < n(g-1)$ then $\MO_X(2n;t)$ is contained in the closure of $\MO_X(2n;t+4)$ in the relevant component $\MO_X(2n)^{\pm}$.
\item For each positive even integer  $t < n(g-1)$, the space $M(V)$ is a single point for a general $V \in \MS_X(2n;t)$. 
\end{enumerate} \end{theorem}

We also compute the dimension of $M(V)$ for general bundles in each component $\MO_{X}(2n)^{\pm}$. It turns out to depend on the class of $n(g-1)$ modulo 4; the precise statement is set out in \S 5.3. A more detailed description of the top strata of $\MO_X(2n)$ is given in  \S 5.4.

To prove these statements, we consider families of symplectic and orthogonal extensions. The main strategy is to relate the invariant $t(V)$ to the geometry of certain higher secant varieties in the projectivized extension spaces. This idea goes back to the work of Lange and Narasimhan \cite{LN}, where Theorem \ref{mainsymp} is proven for rank two bundles using higher secant varieties of the curve $X$ embedded in the extension spaces.

In \S 2, we generalize the geometric framework in \cite{LN} to our situation in Criterion \ref{geomlift}. The embedded curve $X$ in the extension space is replaced by the quadric bundle $\pp E$ and the Grassmannian bundle $\Gr(2, E)$ for symplectic and orthogonal cases respectively. To work with bundles of rank $2n \ge 4$, one has also to understand the situation when two rank $n$ isotropic subbundles intersect non-transversely. It turns out that this case can also be understood geometrically (Criterion \ref{genlift}). This is a key advance upon the methods in \cite{CH2}, which enables us to argue for symplectic and orthogonal bundles of arbitrary rank.

In \S 3, we construct ``universal extension spaces'' parameterizing all the extensions of fixed type, and show that the rational classifying maps to $\MS_X(2n)$ or $\MO_X(2n)$ are defined on dense subsets. In \S 4 and \S 5, we prove the main results for symplectic and orthogonal bundles respectively, using dimension counts based on the geometric information from the preceding sections. In \S 5, we also observe some interesting properties of certain families of orthogonal bundles, stemming from the richer topological structure of the moduli space.\\

\textbf{Acknowledgements:} The first author was supported by the National Research Foundation of Korea (NRF) grant funded by the Korea government (MEST) (No{.}  2011-A423-0004). Also he would like to thank H\o gskolen i Vestfold, Bakkenteigen for the invitation in June 2011.
The second author gratefully acknowledges the hospitality and generous financial support of Konkuk University, Seoul. Both authors thank Christian Pauly and Olivier Serman for helpful discussions.

\section{Symplectic and orthogonal extensions and lifting criteria}

In this section we discuss symplectic and orthogonal extensions. Most of the results in this section were obtained for symplectic extensions in \cite[\S2]{CH2}. Those results will be restated here for the reader's convenience, and the modified versions for orthogonal extensions will be proven in detail. 
\par
Let $V \to X$ be a vector bundle of rank $2n$ equipped with a nondegenerate bilinear form $\omega \colon V \otimes V \to \Ox$, and let $E \subset V$ be a subbundle. Then there is an exact sequence
\[ 0 \to E^{\perp} \to V \to E^{*} \to 0, \]
where  $E^\perp$ is the orthogonal complement of $E$. If $E$ is isotropic of rank $n$, then $E = E^\perp$ and $V$ defines a class $\delta(V) \in H^1(X, \Homom(E^{*}, E)) \cong H^1(X, E \otimes E)$.  

\begin{criterion} Suppose $E$ is simple. An extension $0 \to E \to V \to E^{*} \to 0$ is induced by a symplectic (resp., orthogonal) structure on $V$ with respect to which $E$ is isotropic if and only if the extension class $\delta(V)$ belongs to the subspace
\[ H^{1}(X,  \Sym^{2}E) \quad \left( \hbox{resp.,} \quad H^{1}(X, \wedge^{2}E ) \right) \]
of $H^{1}(X, E \otimes E)$. \label{fundam} \end{criterion}
\begin{proof} This is due to S.\ Ramanan. A detailed proof for the symplectic case is given in \cite[\S2]{Hit1}, and the proof for the orthogonal case is practically identical. \end{proof}

\subsection{Cohomological criterion for lifting}

Here we recall the notion of a bundle-valued principal part (see Kempf \cite{Kem} for corresponding results on line bundles). A locally free sheaf $W$ on $X$ has the flasque resolution
\[ 0 \to W \to \sRat(W) \to \sPrin(W) \to 0, \] 
where $\sRat(W)$ is the sheaf of rational sections of $W$ and $\sPrin(W)$ the sheaf of $W$-valued principal parts. We denote their groups of global sections by $\Rat(W)$ and $\Prin(W)$ respectively. Taking global sections, we obtain
\begin{equation} \label{cohomseq}
0 \to H^0(X, W) \to \Rat(W) \to \Prin(W) \to H^{1}(X, W) \to 0.
\end{equation}
For a principal part $p \in \Prin(W)$, we write its class in $H^{1}(X, W)$ as $[p]$.
\par
Now consider an extension $0 \to E \to V \to E^* \to 0$, and an elementary transformation $F$ of $E^*$ defined by the sequence
\[ 0 \to F \stackrel{\mu}{\to} E^* \to \tau \to 0 \]
for some torsion sheaf $\tau$. We say that $F$ \textsl{lifts to} $V$ if there is a sheaf injection $F \to V$ such that the composition $F \to V \to E^*$ coincides with $\mu$. The following statements are proven in \cite[\S3]{Hit1}:
\begin{lemma} \label{cohomlift} Suppose $h^{0}(X, \Homom(E^{*}, E)) = 0$, and consider an extension $0 \to E \to V \to E^{*} \to 0$ with class $\delta(V) \in H^{1}(X, E \otimes E)$.
\begin{enumerate}
\renewcommand{\labelenumi}{(\arabic{enumi})}
\item There is a bijection between principal parts $p \in \Prin(E \otimes E)$ such that $\delta(V) = [p]$, and elementary transformations of $E^*$ lifting to subbundles of $V$. The bijection is given by
\[ p \ \leftrightarrow \ \Ker \left(p \colon E^{*} \to \sPrin(E) \right). \]
\item Suppose that $\delta(V)  = [p] \in H^{1}(X, \Sym^{2}E)$, corresponding to a symplectic extension. The subbundle corresponding to $\Ker(p)$ is isotropic in $V$ if and only if $p$ is a symmetric principal part; that is, $^tp = p$.
\item Suppose that $\delta(V) = [p] \in H^{1}(X, \wedge^{2}E)$, corresponding to an orthogonal extension. The subbundle corresponding to $\Ker(p)$ is isotropic in $V$ if and only if $p$ is an antisymmetric principal part; that is, $^tp = -p$. \qed
\end{enumerate} \end{lemma}

\subsection{Subvarieties of the extension spaces} \label{alternative}

Given any vector bundle $W$, consider the projectivization $\pi \colon \pp W \to X$. Then we have a natural rational map
\[ \pp W \dashrightarrow \pp H^{1}(X, W) \]
 defined as follows (a slightly different description was given in \cite[\S2.3]{CH2}):

Consider the evaluation map \ $ X \times H^{0}(X, K_X \otimes W^{*}) \to \Kx \otimes W^{*} \label{evalKW} $.
Via Serre duality, the dual of this map is identified with
\[ W \otimes T_{X} \to X \times H^{1}(X, W). \]
Projectivizing this map and then composing with the projection  $X \times \pp H^{1}(X, W) \to  \pp H^{1}(X, W)$, we get a map
\[ \phi \colon \pp W \dashrightarrow \pp H^{1}(X, W). \]
On a fibre $W|_x$, this map is identified with the projectivized coboundary map in the sequence
\begin{equation} 0 \to H^0 (X, W) \to H^{0}(X, W(x)) \to \frac{W(x)}{W}|_{x} \to H^{1}(X, W) \to H^1(X, W(x)) \to 0 . \label{phisequence} \end{equation}
Thus for $w \in \pp W$, the image $\phi(w)$ may be realized as the cohomology class of a $W$-valued principal part supported at $x$ with a simple pole along $w$. As discussed in \cite[\S2.3]{CH2}, the rational map $\phi$ is induced by the complete linear system of the line bundle $\pi^{*}\Kx \otimes \mathcal{O}_{\pp W}(1)$ over $\pp W$.

For $W = \Sym^{2} E$ and $ W= \wedge^2 E$ respectively, we have the rational maps
\[ \pp (\Sym^{2} E) \dashrightarrow \pp H^{1}(X, \Sym^{2}E) \quad \hbox{and} \quad \pp (\wedge^2 E )\dashrightarrow \pp H^{1}(X, \wedge^{2} E). \]
Note that both of these are restrictions of the map $\pp (E \otimes E)  \dashrightarrow \pp H^{1}(X, E \otimes E)$.
In the symplectic case, we consider the chain of maps
\[
\pp E \ \hookrightarrow \   \pp (\Sym^{2} E) \dashrightarrow \pp H^{1}(X, \Sym^{2}E) ,
\]
where the first inclusion is given by the Segre embedding $[v] \mapsto [v \otimes v]$. In the orthogonal case, we consider the chain of maps
\[
\Gr(2,E) \  \hookrightarrow \ \pp (\wedge^{2} E) \dashrightarrow \pp H^{1}(X, \wedge^{2} E),
\]
where the Grassmannian bundle $\Gr(2,E)$ is embedded in $\pp (\wedge^{2} E)$ via the Pl\"ucker embedding.  These rational maps are denoted as
\begin{equation}
 \phi_s \colon \pp E \dashrightarrow \pp H^{1}(X, \Sym^{2}E) \quad \hbox{and} \quad \phi_a: \Gr(2,E) \dashrightarrow \pp H^{1}(X, \wedge^{2} E).
\label{rationalmaps} \end{equation}
\begin{lemma}  \ Let $W = E \otimes E$ for a stable bundle $E$ over a curve of genus $g \ge 2$.
\begin{enumerate}
\renewcommand{\labelenumi}{(\arabic{enumi})}
\item The rational map $\phi \colon \pp (W) \dashrightarrow  \pp H^{1}(X, W)$ is base point free if $\mu(E) < - \frac{1}{2}$, and an embedding if $\mu(E) < -1$. In particular, if $\mu(E) < -1$, the rational maps $\phi_s$ and $\phi_a$ in (\ref{rationalmaps}) are embeddings. 
\item If $E$ is general of negative degree, $\phi$ is injective on a general fiber of $\pp (W)$. If we further assume that $g \ge 3$ then $\phi$ separates two general fibers of $\pp (W)$. 
\end{enumerate} \label{embedding} \end{lemma}
\begin{proof}
(1) (A similar version was stated in \cite[Lemma 2.4]{CH2}, unfortunately with a flawed proof.) One can check that $\phi$ is base point free (resp., an embedding) if $h^{0}(X, W(D)) = 0$ for all effective divisors $D$ of degree one (resp., of degree two) on $X$. Since $W = E \otimes E$ is semistable, these vanishing results follow from the assumption that $\mu(E) < - \frac{1}{2}$ (resp., $\mu(E) < -1$). 

(2) Let $L$ and $M$ be line bundles with $\deg L = \deg E$ and $\deg M = 0$, and put $E_0 = L \oplus M^{\oplus (n-1)}$, where $n = \rank (E)$. One can check that for a general choice of $L$ and $M$, the bundle $E_0 \otimes E_0 (x) $ has no sections. Deforming $E_0$ to a general stable bundle $E$, we see that $h^{0}(X, W(x)) = 0$ for a general $x \in X$. A similar argument shows that if $g \ge 3$, then $h^{0}(X, W(x+y)) = 0$ for general  $ x,y \in X$.  
\end{proof}

\subsection{Irreducibility of the space of principal parts}

In this subsection, we prove some technical facts on principal parts. Let $p$ be a principal part with values in $E \otimes E$.

\begin{definition} The \textsl{degree} $\deg(p)$ of $p$ is defined as the length of the torsion sheaf
\[ \Image \left( p \colon E^{*} \to \sPrin(E) \right). \]
The \textsl{support} of $p$ is defined as the support of $\Image(p)$. A symmetric (resp., antisymmetric) principal part $p$ of degree $k$ (resp., $2k$) is called \textsl{general} if it is supported at $k$ distinct points.
\end{definition}

Obviously, the general symmetric (resp., antisymmetric) principal parts of degree $k$ (resp., $2k$) are parameterized by a quasi-projective irreducible variety. We want to confirm that an arbitrary symmetric or antisymmetric principal part can be obtained as a limit of a continuous family of general ones.

\begin{lemma} {\rm (\cite[Lemma 2.6]{CH2}) } Let $p $ be a $\Sym^2 E$-valued principal part of degree $k$, supported at a single point $x$. Then there exists a local frame $e_{1}, \ldots , e_n$ for $E$ in a neighborhood of $x$, in terms of which $p$ is expressed as
\[ p = \sum_{i=1}^{n} \frac{e_{i} \otimes e_{i}}{z^{k_i}} \]
where $z$ is a uniformizer at $x$, and $ k_{1} \geq k_{2} \geq \cdots \geq k_n \ge 0$, and $k =  \sum_{i=1}^{n} k_{i} $. \qed
\label{sym-irr} \end{lemma}

For antisymmetric principal parts, we have the following analogue:

\begin{lemma} Let $p$ be a $\wedge^{2}E$-valued principal part, supported at a single point $x$. Then there exists a local frame $e_{1}, \ldots , e_n$ for $E$ in a neighborhood of $x$, in terms of which $p$ is expressed as
\[ p = \sum_{i=1}^{s} \frac{e_{2i-1} \wedge e_{2i}}{z^{k_i}} \]
where $z$ is a uniformizer at $x$, and $ k_{1} \geq k_{2} \geq \cdots \geq k_s \ge 0$, and $\deg(p) = 2 \left( \sum_{i=1}^{s} k_{i} \right)$.

In particular, any antisymmetric principal part has even degree.
\label{anti-irr} \end{lemma}
\begin{proof} This argument is adapted from the proof of \cite[Lemma 2.6]{CH2}. Locally, $p$ can be expressed as
\[ p = \frac{1}{z^{k_1}} A \]
for some $n \times n$ antisymmetric matrix $A$ with entries in the ring $R = \cc [z]/ (z^{k_1})$, since we are concerned with the principal parts at $x$ only. In this context, it suffices to show that there exists a matrix $P \in M_{n} (R)$, such that $\det P$ is a unit in $R$ and
\[
^t P A P = \text{diag}(z^{d_1}J, z^{d_2}J, \ldots, z^{d_s}J),
\]
where $0 = d_1 \le d_2 \le \cdots \le d_s$ and $J = \begin{pmatrix} 0 & 1 \\ -1 & 0 \end{pmatrix}$.  This can be shown by mimicking the standard process to get the normal form of an antisymmetric matrix over $\cc$. The only difference lies in that $R$ has non-units contained in the principal ideal $(z)$, and this can be taken care of by allowing the terms $z^{d_1}, \ldots , z^{d_s}$ on the diagonal.

Alternatively, one may work instead over the formal power series ring  $\cc[[z]]$, which is a PID, and truncate the terms in the ideal $(z^{k_1})$ at the final step.  For the process over a PID, see Adkins--Weintraub \cite[Ch.\ 6, Corollary (2.36)]{AW}.

Hence we have a frame for $E$ on this neighborhood in terms of which $p$ appears as
\begin{equation} \sum_{i=1}^{s} \frac{e_{2i-1} \wedge e_{2i}}{z^{k_i}}, \label{antisymppt} \end{equation}
where $k_{i} = k_1 - d_i$. Since the terms in the sum (\ref{antisymppt}) impose independent conditions on sections of $E$, we have
\[ \deg ( p ) = \sum_{i=1}^{s} \deg \left( \frac{e_{2i-1} \wedge e_{2i}}{z^{k_i}} \right) = 2  \sum_{i=1}^{s} k_{i}  , \]
as required. \end{proof}

\begin{corollary} For a fixed vector bundle $E$ and for each $k >0$, the spaces of symmetric principal parts of degree $k$ and antisymmetric principal parts of degree $2k$ are irreducible.
\label{irreducible}  \end{corollary}

\begin{proof} We consider the antisymmetric case; the symmetric case can be proven similarly. It suffices to show the irreducibility of the space of antisymmetric principal parts supported at a single point $x$. By Lemma \ref{anti-irr}, we may choose a trivialization of $E$ near $x$ with respect to which $p$ is expressed as $\sum_{i=1}^{s} p_i$, where
\[ p_{i} = \frac{e_{2i-1} \wedge e_{2i}}{z^{k_i}} \]
and $2 \sum_{i=1}^{s} k_{i} = k$. For each $i$, choose distinct $\lambda^{i}_{1}, \ldots, \lambda^{i}_{k_i} \in \cc$, and define a  family of principal parts
\[ p_{i}(t) = \frac{e_{2i-1} \wedge e_{2i}}{(z - \lambda^{i}_{1}t) \cdots (z - \lambda^{i}_{k_i}t)} \]
where $t$ is a complex parameter. We then put $p(t) = \sum_{i=1}^{s}p_{i}(t)$. By construction, $p(0) = p$, while for small $t \neq 0$ we can rewrite $p(t)$ as a sum of $k$ antisymmetric principal parts of degree 2, supported at $k$ distinct points. Hence we have shown that any antisymmetric principal part of degree $2k$ is a limit of a continuous family of \emph{general} antisymmetric principal parts of degree $2k$.  But we know that the general antisymmetric principal parts are parameterized by an irreducible variety.
This proves the claim.
\end{proof}

\subsection{Geometric criteria for lifting}

In this subsection, we find a geometric interpretation of the cohomological criterion on isotropic liftings in Lemma \ref{cohomlift}. 

Throughout this subsection, $E$ is a general stable bundle of negative degree. In particular, $E$ is simple and $h^0(X, \Homom(E^*, E)) = 0$. Consider the rational maps
\[
\phi_s: \pp E \dashrightarrow \pp H^1(X, \Sym^2 E) \quad \hbox{and} \quad \phi_a: \Gr(2, E) \dashrightarrow \pp H^1(X, \wedge^2 E).
\]
As discussed in Lemma \ref{embedding}, these maps are defined on dense subsets, and furthermore are embeddings if $\mu(E)<-1$. In general, abusing notation, we denote by $\pp E$ and $\Gr(2, E)$ the closures of the images $\phi_s (\pp E)$ in $\pp H^1(X, \Sym^2 E)$ and $\phi_a (\Gr(2, E))$ in $\pp H^1(X, \wedge^2 E)$ respectively.

Next, for a quasi-projective variety $Z \subset \pp^N$, we write $\Sec^k Z$ for the $k$-th secant variety of $\bar{Z}$, which is the closure of the union of linear subspaces spanned by $k$ general points of $Z$.  In particular, $\Sec^1 Z = \bar{Z}$.
\begin{criterion} \label{geomlift} Consider an extension given by
\[ \delta(V) : \quad 0 \to E \to V \to E^{*} \to 0. \]
\begin{enumerate}
\renewcommand{\labelenumi}{(\arabic{enumi})}
\item When $\delta(V) \in \pp H^1(X, \Sym^2 E)$, there is an elementary transformation $F$ of $E^*$ with $\deg (E^*/F) \leq k$ lifting to a Lagrangian subbundle of $V$ if and only if $\delta(V) \in \Sec^{k} \pp E$.\\
\item When $\delta(V) \in \pp H^1(X, \wedge^2 E)$, there is an elementary transformation $F$ of $E^*$ with $\deg (E^*/F) \leq 2k$ lifting to a rank $n$ isotropic subbundle of $V$ if and only if $ \delta(V) \in \Sec^{k} \Gr(2,E)$.  In this case, $\deg E $ and $\deg F$ have the same parity.
\end{enumerate}
\end{criterion}

\begin{proof} In the symplectic case, this is the content of \cite[Lemma 2.10 (2)]{CH2}. The same idea works for antisymmetric case as follows:

By Lemma \ref{cohomlift} (3), an elementary transformation $F \subseteq E^*$ with $\deg (E^*/F) \leq 2k$ lifts to an isotropic subbundle of $V$ if and only if the extension class $\delta(V)$ is of the form $[p]$ where $p$ is an antisymmetric principal part of degree $2l \le k$.

If $p$ is general, it is of the form
\begin{equation} \sum_{i=1}^{l} \frac{e_{i} \wedge f_i}{z_i}, \label{genantisymmppt}
\end{equation}
where $z_1, z_2, \ldots, z_l$ are uniformizers at $l$ distinct points $x_1, x_2,\ldots, x_l \in X$. According to the description of the map $\phi_a$ in \S\ref{alternative}, the class $[p] \in \pp H^1(X, \wedge^2 E)$ lies on the secant plane spanned by $l$ distinct points of $\Gr(2,E)$. Hence $[p] \in \Sec^l \Gr(2, E)$. Conversely, a general point of $\Sec^l \Gr(2, E)$ corresponds to a class $\delta(V) =[p]$ where $p$ is of the form (\ref{genantisymmppt}).

From this correspondence on the general points, we get the desired statement by 
Corollary \ref{irreducible}. In this case, $F = \Ker \left(p \colon E^{*} \to \sPrin(E) \right)$ with $\deg p = 2l$, so $\deg F$ has the same parity as $\deg E$.
\end{proof}
\begin{corollary}
Consider an extension given by
\[ \delta(V) : \quad 0 \to E \to V \to E^{*} \to 0. \]
\begin{enumerate}
\renewcommand{\labelenumi}{(\arabic{enumi})}
\item Assume $\delta(V) \in H^1(X, \Sym^2 E)$. If $\delta(V) \le \Sec^k \pp E$, then $t(V) \le 2(k + \deg E)$.\\
\item Assume $\delta(V) \in H^1(X, \wedge^2 E)$. If $\delta(V) \le \Sec^k \Gr(2, E)$, then $t(V) \le 2(2k + \deg E)$.
\end{enumerate}
\end{corollary}
\begin{proof}
By definition of the invariant, $t(V) \le - 2 \deg F$ for any rank $n$ isotropic subbundle $F \subset V$. The above bounds follow as a direct consequence of Criterion \ref{geomlift}.
\end{proof}
\par
The converse of the above Corollary is not true in general. For instance, assume $\delta(V) \in H^1(X, \Sym^2 E)$. The bound $t(V) \le 2(k + \deg E)$ tells us that there is a Lagrangian subbundle $F$ of degree $\ge -\deg E -k$, but this does not imply that $F$ lifts from an elementary transformation of $E^*$. So in general, we are led to consider a diagram of the form
\begin{equation}
\xymatrix{ 0 \ar[r] & E \ar[r] & V \ar[r] & E^{*} \ar[r] & 0 \\
 0 \ar[r] & H \ar[r] \ar[u] & F \ar[r] \ar[u] & G \ar[r] \ar[u] & 0 } \label{nontransv}
 \end{equation}
where $H \subset E$ is a subbundle of degree $-h \leq 0$ and rank $r \geq 0$, and $G$ is a locally free subsheaf of $E^*$ of rank $n-r$. When $r = 0$ so that $E|_x$ and $F|_x$ meet transversely for general $x \in X$, this reduces to the situation of Criterion \ref{geomlift} of the lifting of elementary transformations. The remaining part of this subsection will be devoted to finding a criterion for the existence of such a diagram with $r>0$.

For $H \subseteq E$, let $q \colon E \to E/H$ be the quotient map.  Since $H$ is isotropic, $H^{\perp}$ fits into the diagram
\[
\xymatrix{ 0 \ar[r] & E \ar[r] & V \ar[r] & E^{*} \ar[r] & 0 \\
 0 \ar[r] & E \ar[r] \ar@{=}[u] & H^{\perp} \ar[r] \ar[u] & \left( E/H \right)^{*} \ar[r] \ar[u]^{^tq} & 0. }
 \]
For $\delta (V) \in H^1(X, \Homom(E^*, E))$, we have
\[
\delta \left( H^{\perp} \right) = \, ^tq^{*} (\delta(V)) \ \in \ H^{1}(X, \Homom((E/H)^{*}, E)).
\]
Furthermore, $H^{\perp}$ inherits a (degenerate) bilinear form from $V$. Since $(H^{\perp})^{\perp} = H$, the quotient $H^{\perp}/H$ is equipped with a \emph{nondegenerate} bilinear form coming from $V$. In fact, $H^{\perp}/H$ is a symplectic (resp., orthogonal) extension in the upper exact sequence of
\[
\xymatrix{ 0 \ar[r] & E/H \ar[r] & H^{\perp} /H  \ar[r] & \left( E/H \right)^{*} \ar[r] & 0 \\
 0 \ar[r] & E \ar[r] \ar[u]^q & H^{\perp} \ar[r] \ar[u] & \left( E/H \right)^{*} \ar[r] \ar@{=}[u] & 0 }
 \]
corresponding to the class
\[
\delta(H^{\perp}/H) = q_{*}(^tq)^{*} (\delta(V))
\]
in $ H^{1} \left( X, \Sym^{2} (E/H) \right)$  (\text{resp.,} $H^{1} \left( X, \wedge^{2} (E/H) \right)$ ).

\begin{lemma} The induced maps
\[ q_{*}(^tq)^{*} \colon H^{1}(X, \Sym^{2} E ) \to H^{1} \left( X, \Sym^{2} (E/H) \right) \]
and
\[ q_{*}(^tq)^{*} \colon H^{1}(X, \wedge^{2} E ) \to H^{1} \left( X, \wedge^{2} (E/H) \right) \]
are surjective.
\label{qq} \end{lemma}
\begin{proof} These maps are induced from
\[ (^tq)^{*} : E \otimes E \to (E/H) \otimes E \quad \hbox{and} \quad q_{*} \colon (E/H) \otimes E \to (E/H) \otimes (E/H) \]
respectively. By local computation, it can be seen that the images of $q_{*}(^tq)^{*}$ in $(E/H) \otimes (E/H)$ are precisely $\Sym^{2} \left( E/H \right)$ and $\wedge^{2} \left( E/H \right)$ respectively. Hence on the first cohomology level, the induced maps are surjective. \end{proof}

We obtain the following geometric criterion on liftings:
\begin{criterion} \label{genlift} Let $V$ be an extension of $E^*$ by $E$ with class $\delta (V)$.  Fix a subbundle $0 \neq H \subset E$ and write $\deg(E) = -e, \deg(H) = -h$. Let $f > 2h - e$.
\begin{enumerate}
\renewcommand{\labelenumi}{(\arabic{enumi})}
\item Assume $\delta(V) \in H^1(X, \Sym^2 E)$. Then $V$ admits a Lagrangian subbundle $F$ of degree $\ge -f$ inducing a diagram of the form (\ref{nontransv}) if and only if
\[
q_{*}(^tq)^{*} (\delta(V)) \ \in \   \Sec^{(e+f -2h)} \pp(E/H) .
\]
\item Suppose $\delta(V) \in H^1(X, \wedge^2 E)$ and $e \equiv f \mod 2$. Then $V$ admits a rank $n$ isotropic subbundle $F$ of degree $\ge -f$ inducing a diagram of the form (\ref{nontransv}) if and only if
\[
q_{*}(^tq)^{*} (\delta(V)) \ \in \   \Sec^{\frac{1}{2}(e+f - 2h)} \Gr(2, E/H) .
\]
\end{enumerate}
\end{criterion}

\begin{proof} Since the arguments for (1) and (2) are parallel, we prove (2) only. The orthogonal bundle $V$ admits an isotropic subbundle $F$ inducing the diagram (\ref{nontransv}) if and only if $H^{\perp}$ admits a subbundle $F$ of degree $\ge -f$ which is isotropic with respect to the antisymmetric form inherited from $V$, yielding the diagram
\[
\xymatrix{ 0 \ar[r] & E \ar[r] & H^{\perp}  \ar[r] & \left( E/H \right)^{*} \ar[r] & 0 \\
 0 \ar[r] & H \ar[r] \ar[u] & F \ar[r] \ar[u] & G \ar[r] \ar[u] & 0 }
 \]
Since $G = F/H$ has the same rank as $(E/H)^*$, the map $G \to (E/H)^*$ is an elementary transformation whose quotient is a torsion sheaf of degree $\le e+f - 2h$.  Factorizing by $H$, we get
\[
\xymatrix{ 0 \ar[r] & E/H \ar[r] & H^{\perp}/H   \ar[r] & \left( E/H \right)^{*} \ar[r] & 0 \\
  &   & F/H \ar[r]^= \ar[u] & G \ar[r] \ar[u] & 0. }
 \]
We are in this situation precisely when the orthogonal extension $H^{\perp}/H$ with class
\[ q_{*}(^tq)^{*} \delta(V) \in H^1 (X, \wedge^2 (E/H)) \]
admits an isotropic lifting of some elementary transformation of $(E/H)^*$ of degree $\le e+f - 2h$. By Criterion \ref{geomlift} (2), this is equivalent to the class $\delta \left( H^{\perp}/H \right)$ belonging to $\Sec^{\frac{1}{2}(e+f -2h)} \Gr(2,E/H)$ in $\pp H^{1}(X, \wedge^{2}(E/H))$.
\end{proof}

\section{Parameter spaces of extensions}

Here we construct ``universal extension spaces'', following Lange \cite{Lan2} (see also \cite[\S4.2]{CH2}), and investigate stability of the corresponding symplectic and orthogonal bundles.

\subsection{Construction of the families}

For a positive integer $e$, let $U_X(n,-e)^s$ denote the moduli space of stable vector bundles of rank $n$ and degree $-e < 0$ over $X$. By Narasimhan--Ramanan \cite[Proposition 2.4]{NR1975}, there exist a finite \'etale cover
\[ \pi_{e} \colon \tU_e \to U_X(n, -e)^s \]
and a bundle $\cE_e \to \tU_e \times X$ with the property that $\cE_e|_{\{ E \} \times X} \cong \pi_e(E)$ for all $E \in \tU_e$. When $\gcd(n, e)=1$, it is well known that $\pi_e$ reduces to the identity map.

By Riemann-Roch and semistability, for each $E \in U_X(n, -e)^s$, we have
\[ h^1(X, \Sym^{2}E) = (n+1)e + \frac{n(n+1)}{2}(g - 1). \]
Therefore, the sheaf $R^{1} {p}_{*} (\Sym^2 (\cE_e))$ is locally free of rank $(n+1)(e + \frac{1}{2}n(g - 1))$ on $\tU_e$. We denote its projectivization by $\mu \colon \sy_{e} \to \tU_e$. We have a diagram
\begin{displaymath}
\label{universal} \xymatrix{
\sy_{e}  \ar[dd]_{\mu} & \sy_{e} \times X \ar[d]^{\mu \times \Iden_X} \ar[l]_{r} & \\
& \tU_{e} \times X \ar[dl]_{p} \ar[dr]^{q} & \\
\tU_{e} & & X}
\end{displaymath}
We write $r \colon \sy_e \times X \to \sy_e$ for the projection. By Lange \cite[Corollary 4.5]{Lan2}, there is an exact sequence of vector
bundles
\[ 0 \to (\mu \times \Iden_{X})^{*}\cE_e \otimes r^{*}\mathcal{O}_{\sy_e}(1) \to \cW_{e} \to (\mu \times \Iden_{X})^{*}\cE_e^{*} \to 0 \]
over $\sy_{e} \times X$, with the property that for $\delta \in \sy_e$ with $\mu(\delta) = E$, the restriction of $\cV_{e}$ to $\{ \delta \} \times X$ is isomorphic to the extension of $E^*$ by $E$ defined by $\delta \in \pp H^1(X, \Sym^{2}E)$. By Lemma \ref{cohomlift}, the space $\sy_{e} $ classifies all symplectic extensions of $E^*$ by $E$ for all $E \in U(n, -e)^s$, up to homothety.

In the same way, we define a bundle $\A_{e} \to \tU_e$ whose fibre at $E$ is $\pp H^{1}(X, \wedge^{2}E)$, a projective space of dimension $(n-1)(e + \frac{1}{2}n(g - 1)) - 1$. There is a sequence of vector bundles
\[ 0 \to (\mu \times \Iden_{X})^{*}\cE_e \otimes r^{*}\mathcal{O}_{\A_e}(1) \to \cV_{e} \to (\mu \times \Iden_{X})^{*}\cE_e^{*} \to 0 \]
over $\A_{e} \times X$, with the property that for $\delta \in \A_e$ with $\mu(\delta) = E$, the restriction of $\cW_{e}$ to $\{ \delta \} \times X$ is isomorphic to the extension of $E^*$ by $E$ defined by $\delta \in \pp H^1(X, \wedge^{2}E)$. Again by Lemma \ref{cohomlift}, the space $\A_{e} $ classifies all the orthogonal extensions of $E^*$ by $E$ for all $E \in U(n, -e)^s$, up to homothety.

For each $e>0$, the universal bundles $\cW_e \to \sy_e \times X$ and $\cV_e \to \A_e \times X$ induce classifying maps
\[
\sigma_e: \sy_e \dashrightarrow \MS_X(2n)
\ \ \text{and} \ \
\alpha_e: \A_e \dashrightarrow \MO_X(2n)
\]
respectively. The indeterminacy loci of these maps consist of precisely the points whose associated symplectic/orthogonal bundles are not semistable.

\subsection{Stability of extensions}

The universal extension spaces $\sy_e$ and $\A_e$ provide a natural way to study the Segre stratification on the moduli space of symplectic/orthogonal bundles, via the classifying maps $\sigma_e$ and $\alpha_e$. In order to proceed in this direction, we must verify that a general bundle with extension class represented in $\sy_e$ (resp., $\A_e$) is a stable symplectic (resp., stable orthogonal) bundle.

The same question for vector bundles was formulated by Lange \cite{Lan}, and solved in Brambila-Paz--Lange \cite{BPL} and Russo--Teixidor i Bigas \cite{RT}. In both papers, elementary transformations were used to construct stable bundles with the prescribed Segre invariant. In this subsection, we prove the existence of a stable orthogonal/symplectic bundle in $\sy_e$ and $\A_e$ for each $e >0$. Elementary transformations are used, but in a somewhat different context.

We begin by establishing the statement for $e = 1$ and $e = 2$. We will need the following bound on the classical Segre invariants of vector bundles, due to Hirschowitz \cite[Th\'eor\`eme 4.4]{Hir}:

\begin{lemma} \label{maxsegre} Let $E$ be a general stable vector bundle of rank $n$ and degree $-e$. Let $H \subset E$ be a subbundle of rank $r$ and degree $-h$. Then
\[
r(n-r) (g-1) \ \le \ nh -re \ < \ r(n-r) (g-1) +n .
\]
\end{lemma}
\noindent Also we need the following result called the Hirschowitz lemma and its variant.

\begin{lemma}  \label{Hirlemma}   \ \ 
\begin{enumerate}
\renewcommand{\labelenumi}{(\arabic{enumi})}
\item Let $H_1$ and $H_2$ be general stable bundles, $\rank (H_i) = r_i$  and $\deg (H_i) = d_i$ for $i = 1,2$. If $r_1 d_2 + r_2 d_1 \ge r_1 r_2(g-1)$ so that $\mu(H_1 \otimes H_2)  \ge g-1$, then $h^1(X, H_1 \otimes H_2) = 0$.
\item Let $F$ be a general stable bundle of rank $n$. If $\deg F \ge \frac{1}{2}n(g-1)$ so that $\mu(F \otimes F) \ge g-1$, then $h^1(X, F \otimes F) = 0$. \end{enumerate} \end{lemma}
\begin{proof}
Part (1) was proven by Hirschowitz \cite[4.6]{Hir}; see also Russo--Teixidor i Bigas \cite[Theorem 1.2]{RT}. The variant (2) is \cite[Lemma A1]{CH2}.
\end{proof}

\begin{proposition} \ Let $E \in U(n, -e)^s$ be general, $e = 1,2$.
\begin{enumerate}
\renewcommand{\labelenumi}{(\arabic{enumi})}
\item For $e=1$, every point in $\pp H^1(X, \Sym^2 E)$  outside a sublocus $Y$ with $\dim Y \le n = \dim \pp E$ corresponds to a stable symplectic bundle. Hence a general point of $\sy_1$ represents a stable symplectic bundle.
\item For $e=1$, every point in $\pp H^1(X, \wedge^2 E)$ corresponds to a stable orthogonal bundle. For $e=2$, every point of $\pp H^1(X, \wedge^2 E)$ outside a sublocus $Z$ with $\dim Z \le 2(n-2) +1 = \dim Gr(2,E)$ corresponds to a stable orthogonal bundle. Hence a general point of 
 $\A_1$ or $\A_2$ represents a stable orthogonal bundle.
\end{enumerate}
\label{stabilityonetwo} \end{proposition}

\begin{proof} Consider a symplectic or orthogonal extension $0 \to E \to V \to E^{*} \to 0$ where $E$ is a general stable bundle of rank $n$ and degree $-e \in \{-1,-2 \}$.  
Assume that $V$ is not stable, so there is  an isotropic subbundle $F$ of $V$ of rank $r (\leq n)$ and degree $\ge 0$. 
The intersection of $E$ and $F$ contains a subbundle $F_1$ of rank $r_1$ (possibly zero), and the image of $F$ in $E^*$ is a locally free subsheaf $F_2$ of rank $r_2$, yielding a diagram
\begin{equation} \label{notstable} \xymatrix{ 0 \ar[r] & E \ar[r] & V \ar[r] & E^{*} \ar[r] & 0 \\
0 \ar[r] & F_{1} \ar[r] \ar[u] & F \ar[r] \ar[u] & F_{2} \ar[r] \ar[u] & 0. } \end{equation}
We will bound the dimension of the locus of the extensions $V$ admitting this kind of diagram. 

First assume $r_{1} \neq 0$ and so $r_2 < n$. Since $E$ is general, by Lemma \ref{maxsegre} we have
\begin{equation} \label{degF_1} \deg F_{1} \leq - \frac{r_1}{n} \left( e + (n-r_{1})(g-1) \right)  \ <0 .
\end{equation}
In the same way (if $r_2 \neq 0$),
\begin{equation} \label{degF_2}
 \deg F_{2} \leq \frac{r_2}{n} \left( e - (n-r_{2})(g-1) \right),
 \end{equation}
which implies $\deg F_2 \le 0$ since $r_{2} < n$. Therefore, $\deg F < 0$.

Next, assume $r_1 = 0$ and $ r = r_{2}  = n$.
In this case, $F$ is an elementary transformation of $E^*$. If $\deg F \ge 0$, then the torsion sheaf $E/F$ has degree $\le e$.
In the symplectic case, we need only consider the case  when $e = 1$. By Criterion \ref{geomlift} (1), if $F$ lifts to $V$ as a Lagrangian subbundle, then $\delta(V) \in \phi_s(\pp E ) $ in $\pp H^1(X, \Sym^2 E)$. 
In the orthogonal case, $e \le 2$. By Criterion \ref{geomlift} (2), if $F$ lifts to $V$ as an isotropic subbundle, then in fact $e=2$ and $\delta(V) \in { \Gr(2, E) }$ in $\pp H^1(X, \wedge^2 E)$.  

Finally assume $r_1 = 0$ and $r=r_2 < n$. From the inequality (\ref{degF_2}),  $\deg F = \deg F_2 \le 0$. The possibility $\deg F = 0$ appears only for the following special cases:
\begin{enumerate}
\renewcommand{\labelenumi}{(\roman{enumi})}
\item $e = 1$; $r = n - 1$; $g = 2$
\item $e = 2$; $r = n - 1$; $g = 3$
\item $e = 2$; $r = n - 1$; $g = 2$
\item $e = 2$; $r = n - 2$; $g = 2$.
\end{enumerate}
From now on, we show that  in each of these cases, the dimension of the locus of extensions admitting a lifting of a subsheaf $F$ of $E^*$ as an isotropic subbundle of $V$ of degree zero and rank $r<n$  is bounded by $n= \dim \pp E$ in $\pp H^1(X, \Sym^2 E)$ and $2(n-2) +1 = \dim Gr(2,E)$ in $\pp H^1 (X, \wedge^2 E)$ respectively.

Note that in each case (i)--(iv), $F$ is a maximal subbundle of $E^*$, which is easily seen from the inequalities in Lemma \ref{maxsegre}. Hence the quotient $Q := E^{*}/F$ is torsion-free of degree $e$ and we obtain the diagrams
\begin{equation}
\xymatrix{ 0 \ar[r] & E \ar[r] & \left( F^{\perp} \right)^{*} \ar[r] & Q \ar[r] & 0 \\
0 \ar[r] & E \ar[r] \ar[u]^{=} & V \ar[r] \ar[u] & E^{*} \ar[r] \ar[u] & 0 \\
 & & F \ar[r]^{=} \ar[u] & F \ar[u] & }
\label{V} \end{equation}
and
\begin{equation}
 \xymatrix{ 0 \ar[r] & Q^{*} \ar[r] & G \ar[r] & Q \ar[r] & 0 \\
0 \ar[r] & Q^{*} \ar[r] \ar[u]^{=} & F^{\perp} \ar[r] \ar[u] & E^{*} \ar[r] \ar[u] & 0 \\
 & & F \ar[r]^{=} \ar[u] & F. \ar[u] & }
 \label{Fperp} \end{equation}
Since $(F^{\perp})^{\perp} = F$, the bundle $G = F^{\perp}/F$ inherits the nondegenerate bilinear form from $V$. Since $F^{\perp} \cap E = Q^*$ is contained in the isotropic subbundle $E \subset V$, it is an isotropic subbundle of $G$. Thus the class of the extension $G$ belongs to either $H^1(X, \Sym^2 Q^*)$ or $H^{1}(X, \wedge^{2} Q^{*})$ by Criterion \ref{fundam}.
\par
For a given general $E$, there are finitely many choices for $F$ in cases (i), (ii), and (iv), while in case (iii), the subbundles of degree zero and rank $n - 1$ in $E$ vary in a Quot scheme of dimension $n - 1$. Once $F$ is chosen, the quotient $Q = E^{*} / F$ is fixed. After we choose a symplectic or orthogonal extension $G$ of $Q$ by $Q^*$, the bundles $F^\perp$ and $V$ are determined from the above diagrams by
\[
(F^{\perp})^* = \left( E \oplus G \right) / Q^*
\]
and
\[
 V = \left( E \oplus F^{\perp} \right) / Q^{*}.
 \]
Therefore, for a fixed general $E$, the dimension of the locus of $V$ appearing in the above class of diagram is bounded by that of the deformations of $F$ and $G$.

Let us consider the orthogonal case first. In cases (i)--(iii), the bundle $Q$ has rank 1, so there is no nontrivial orthogonal extension of $Q$ by $Q^*$. Hence the only possibility for $V$ is the direct sum $E \oplus E^*$, which is excluded.
In case (iv), we have $h^1(X, \wedge^2 Q^*) = 3$, and so
\[ \dim \{ F : F \subset E^* \} + h^{1}(X, \wedge^{2}Q^{*}) - 1 = 2. \]
Since $0 <r = n-2$, we have $2 \le 2(n-2)+1$ as was claimed.
\par
For the symplectic case, we assumed $e = 1$, so (i) is the only case to be considered. In this case, there are finitely many choices of $F$, and $h^{1}(X, \Sym^{2}Q^{*}) - 1 = 1 \le n$, as was claimed.

This confirms that a general point of $\sy_1$, $\A_1$ or $\A_2$ represents a stable orthogonal bundle.
\end{proof}
 
\begin{theorem} For each $e > 0$, a general point of $\sy_e$ (resp., $\A_e$) represents a stable symplectic (resp., orthogonal) bundle. \label{stability} \end{theorem}
\begin{proof} We consider the orthogonal case first: For each value of $e>0$, we will exhibit a stable orthogonal bundle represented in $\A_e$. The statement will then follow from the openness of the stable objects in families.

Let $E \in U_{X}(n,-1)$ be general.
For any $k \geq 2$, choose a general antisymmetric principal part $p$ of degree $2k$, which defines an element of $\Sec^k \Gr(2, E)$. Since $\Gr(2,E)$ is nondegenerate and properly contained in $\Sec^{k} \Gr(2,E)$, we may assume that $[p]$ does not lie on the image of $\Gr(2, E)$. By Lemma \ref{cohomlift}, the sheaf $F := \Ker \left( p \colon E^{*} \to \sPrin(E) \right)$ lifts to a rank $n$ isotropic subbundle of the extension $V$ with $\delta(V) = [p]$. Deforming $p$ if necessary, we can assume that $F$ is stable. By Proposition \ref{stabilityonetwo} (2), moreover, $V$ is a stable orthogonal bundle. Since $V$ fits into an orthogonal extension
\begin{equation} 0 \to F \to V \to F^{*} \to 0, \label{extF} \end{equation}
it is associated to a point in $\A_{2k-1}$.

In the same way, consider a general bundle $E \in U_X(n, -2)$ and choose a general antisymmetric principal part of degree $2k+2$ for each $k \geq 2$, defining a point of $\Sec^{k+1} \Gr(2,E)$. By Proposition \ref{stabilityonetwo} (2), we may assume that the extension $V$ associated to the point $[p]$ is a stable orthogonal bundle, if we avoid a subvariety $Z$ with $\dim Z \le \dim \Gr(2,E)$; and this is possible since $\Gr(2,E)$ is nondegenerate and properly contained in $\Sec^{k+1} \Gr(2,E)$. As in the previous case, $V$ fits into an orthogonal extension (\ref{extF}) where this time $\deg(F) = -2k$, so $V$ is represented in $\A_{2k}$.
\par
In the symplectic case, we argue similarly. For each $k \geq 2$, we can choose a general symmetric principal part $q$ of degree $k+1$ which defines a point of $\Sec^{k+1} \pp E $. By the above argument, the extension $V$ determined by $q$ is a stable symplectic bundle which is represented in $\sy_{k}$. \end{proof}

\section{Description of the Segre strata for symplectic bundles}

In this section, we use the map $\sigma_e: \sy_e \dashrightarrow \MS_X(2n)$ discussed in the previous section to prove Theorem \ref{mainsymp}.
Recall that for each positive even integer $t$, we denote by $\MS_X(2n;t)$ the sublocus of $\MS_X(2n)$ consisting of symplectic bundles $V$ with $t(V) = t$. Also, for $V \in \MS_X(2n)$, we write $M(V)$ for the space of Lagrangian subbundles of $V$ of (maximal) degree $-\frac{1}{2}t(V)$.

\begin{theorem} \label{dimension} \ Consider a positive even integer $ t=2e \le  n (g-1)+1$.
\begin{enumerate}
\renewcommand{\labelenumi}{(\arabic{enumi})}
\item A general point of $\sy_e$ corresponds to a bundle $V$ with $t(V) = t$; in particular, $\MS_X(2n;t)$ is nonempty. Furthermore, $\MS_X(2n;t)$ is irreducible.
\item If $t < n(g-1)$, then $M(V)$ is a single point for a general $V \in \MS_X(2n;t)$.
\item If $t \in \{n(g-1), \: n(g-1) +1 \}$ and $V$ is general in $\MS_X(2n;t)$ then any pair of Lagrangian subbundles in $M(V)$ intersect transversely in $V$ in a general fiber.
\end{enumerate}
\end{theorem}
\begin{remark}
Statements (2) and (3) can be viewed as a symplectic analogue of Lange--Newstead \cite[Theorem 2.3 \& Proposition 2.4]{LaNe}.
\end{remark}
\begin{proof} Let $E$ be a general bundle in $U(n, -e)^s$ for $e \ge 1$. By Theorem \ref{stability}, the moduli map $\sy_{e}|_{E} \cong \pp H^{1}(X, \Sym^{2}E) \dashrightarrow \MS_X(2n)$ is defined on a dense subset.
We want to compute a bound on the dimension of the locus in $H^{1}(X, \Sym^{2}E)$ of extensions $V$ which admit a Lagrangian subbundle $F$ of degree $\ge -e$ other than $E$. There are two possibilities: either $F$ is an elementary transformation of $E^*$ lifting to $V$, or $F$ fits into a diagram of the form (\ref{nontransv}).

\textit{Step 1.} We show first that the latter situation does not arise for a general $V$ in $H^1(X, \Sym^2 E)$ for $t = 2e \le n(g-1) +1$. As before, we write $\deg(H) = -h$. Recall from the diagram (\ref{nontransv}) that $G$ is an elementary transformation of $(E/H)^*$ whose quotient is a torsion sheaf of degree
\[
\deg (E/H)^* - \deg G = (e-h) - (\deg F + h) \le 2(e-h).
\]
In particular,
\begin{equation}
e-h \ge 0.
\label{eandh}
\end{equation}
Now for each fixed $H \subset E$, Criterion \ref{genlift} says that the locus of extensions in $H^1(X, \Sym^2 E)$ admitting a diagram of the form (\ref{nontransv}) for some $F$ of degree $\ge -e$ is bounded by
\[ \dim \left( \Sec^{2(e-h)}\pp(E/H) \right) + 1 + \dim \Ker \left( q_{*} \, ^tq^{*} \right). \]
The secant variety has dimension bounded by
\[ 2(e-h)(n-r) + 2(e-h) - 1 = 2(e-h)(n-r+1) - 1. \]
Also by Proposition \ref{qq}, we have
\begin{align*} \dim \Ker \left(q_{*} \, ^tq^{*} \right) &= h^{1}(X, \Sym^{2}E) - h^{1} \left(X, \Sym^{2}(E/H) \right) \\
 &\leq h^{1}(X, \Sym^{2}E) - (n-r+1)(e-h) - \frac{1}{2} (n-r)(n-r+1)(g-1). \end{align*}
Finally we take account of the deformations of $H$ by computing the dimension of the appropriate Quot scheme of $E$. Since $E$ is general, by Lemma \ref{maxsegre} we have
\[
 nh - re \ \ge \ r(n-r)(g-1).
\]
We may furthermore assume that $H$ and $E/H$ are general, and so
\[ h^1(X, \Homom (H, E/H)) = 0 \]
by Lemma \ref{Hirlemma} (1). Therefore $[H]$ is a smooth point of the Quot scheme (cf.\ Le Potier \cite[p.\ 125]{LeP}), and the dimension of the Quot scheme is given by
\[ h^{0}(X, \Homom(H, E/H)) = nh - re - r(n-r)(g-1). \]
Adding up these three terms, we see that the dimension of the locus of extensions in $H^1(X, \Sym^2 E)$ admitting a diagram of the form (\ref{nontransv}) for some $F$ of degree $\ge -e$ is bounded by
\begin{eqnarray*}
\dim \left( \Sec^{2(e-h)} \pp (E/H) \right) + 1 + \dim \Ker (q_* \, ^t q^*) + h^0(X, \Homom (H, E/H)) \\
\le \ (e-h)(n-r+1) - \frac{1}{2} (n-r)(n+r+1) (g-1) + nh-re + h^1(X, \Sym^2 E) \\
= \ e(n-r) - (e-h)(r-1) - \frac{1}{2} (n-r)(n+r+1) (g-1) + h^1(X, \Sym^2 E).
\end{eqnarray*}
By the assumption $e \le \frac{1}{2} (n (g-1)+1)$ and the inequality (\ref{eandh}), this is smaller than the dimension of the whole extension space $h^1(X, \Sym^2 E)$. This shows that for $t = 2e \le  n (g-1) +1$, a general $V \in H^1(X, \Sym^2 E)$ does not admit a Lagrangian subbundle $F$ of degree $\ge -e$ which fits into a diagram of the form (\ref{nontransv}).

\textit{Step 2.} Next we consider the case $r = 0$, so $F$ is an elementary transformation of $E^*$ lifting to $V$ isotropically. Let $\deg F = -f \ge -e$. By Criterion \ref{geomlift}, the dimension of the locus of extensions in $H^1(X, \Sym^2 E)$ admitting such an $F$ is bounded by
\[
\dim \left( \Sec^{e+f} \pp E \right) +1 \ \le \ (e+f)(n+1).
\]
If $f < \frac{1}{2} n (g-1)$, this bound is smaller than
\[
 h^1(X, \Sym^2 E)  = (n+1)e + \frac{1}{2}n(n+1)(g-1).
\]

\textit{Step 3.} Combining the dimension counts for two possibilities in the above, we conclude:\\
(i) if $t=2e <  n (g-1)$, a general $V \in H^1(X, \Sym^2 E)$ does not have a Lagrangian subbundle of degree $\ge -e$ other than $E$ itself. This shows (2). \\
(ii)  if $t=2e \in \{ n (g-1), \: n(g-1) + 1 \}$, a general $V \in H^1(X, \Sym^2 E)$ does not have a Lagrangian subbundle of degree $> -e$. Also, any Lagrangian subbundle of degree $-e$ different from $E$ intersects $E$ transversely at a general fiber.  This shows (3).

Hence for a general $V \in H^1(X, \Sym^2 E)$ with $t = 2e \le n(g-1) +1$, we have $t(V) = t$. It follows that the rational map $\sigma_e : \sy_e \dashrightarrow \MS_X(2n)$ sends a general point of $\sy_e$ to $\MS_X(2n;t)$, so $\MS_X(2n;t)$ is non-empty for each $t \le n (g-1) +1$.  Since $\sy_e$ is irreducible, so is its image. To see that $\MS_X(2n;t)$ is irreducible, it suffices to show that every point of $\MS_X(2n;t)$ is in the closure of the image of $\sy_e$. Any $V$ in $\MS_X(2n;t)$ admits a Lagrangian subbundle $E$ of degree $-e$ which might be unstable. But every such $E$ is contained in an irreducible family of bundles whose general member is  a stable bundle in $U_X(n, -e)$. This shows (1).
\end{proof}

To finish the proof of Theorem \ref{mainsymp}, we need the following description of the tangent space of $M(V)$.

\begin{lemma} \label{Zariskitangent}
For a symplectic bundle $V \in \MS_X(2n)$ and a point $[\eta : E \subset V]$ of $M(V)$, the Zariski tangent space of $M(V)$ at $\eta$ is identified with $H^0(X, \Sym^2 E^*)$.
\end{lemma}
\begin{proof}
It is well known that the Zariski tangent space to the Quot scheme of a vector bundle $V$ at a point $[E\subset V]$ is given by $H^0(X, \Homom(E, V/E))$. Intuitively, this can be explained as follows: A rank $n$ subbundle of $V$ is equivalent to a global section of the Grassmannian bundle $\Gr(n, V)$ over $X$. A tangent vector to the Quot scheme of $V$ at $[E \subset V]$ corresponds, at each fiber $x \in X$, to a tangent vector to the Grassmannian $\Gr(n, V_x)$ at $[E_x \subset V_x]$. Therefore, the tangent vector is a global section of the bundle $\Homom(E, V/E)$. When $E$ is a Lagrangian subbundle of a symplectic bundle $V$, we have $V/E \cong E^*$ and so $\Homom(E, V/E) \cong E^* \otimes E^*$. In this case, a similar argument shows that a tangent vector to $M(V)$ at $[E \subset V]$ corresponds to a global section of $\Sym^2 E^*$, since for each $x \in X$ the tangent space of the Lagrangian Grassmannian at $[E_x \subset V_x]$ is identified with $\Sym^{2}E^{*}|_x$.
\end{proof}

\noindent \textit{Proof of Theorem \ref{mainsymp}}\\ 
\\
We begin with part (3) of the statement. The first part of (3) is Theorem \ref{dimension} (2). For the latter part, we invoke Lemma \ref{Zariskitangent}: $\dim M(V) = h^0(X, \Sym^2 E^*)$ for a smooth point $[E \subset V]$ of $M(V)$. By Lemma \ref{Hirlemma} (2), for a general $E \in U_X(n,-e)$, we have
 \[
\dim H^0(X, \Sym^2 E^*) \ = \ \begin{cases}
0&  \text{if} \  2e =  n(g-1),  \\
\frac{n+1}{2}&  \text{if} \  2e=n(g-1)+1.
\end{cases}
\] 
(1)  A straightforward computation   shows that when $  n(g-1) \le 2e \le n(g-1)+1 $,
\[ \dim M(V) = \dim \sy_e - \dim \MS_X(2n). \]
Since  $ \dim M(V)  \ge \dim\sigma_{e}^{-1}(V)$, this equality  implies that $\sigma_e$ is dominant. This shows that for a general $V \in \MS_X(2n)$, we have
\[
n(g-1) \le t(V) \le n(g-1) +1.
\]
By semicontinuity, $t(V) \le n(g-1) +1$  for any symplectic bundle $V$.\\
\\
(2) If $t = 2e \ge n(g-1) $, then $\dim \MS_X(2n;t) = \dim \MS_X(2n)$ by (1). 
Assume $t = 2e < n(g-1)$. By Theorem \ref{dimension} (2), the map $\sigma_e \colon \sy_{e} \dashrightarrow \MS_X(2n;t)$ is generically finite (of degree $\deg \pi_e$), and so
\[
\dim \MS_X(2n;t) \ = \ \dim \sy_e \ = \ \frac{1}{2}n(3n+1)(g-1) + (n+1)e.
\]

For the last part of (2), we can use Criterion \ref{geomlift} to show that the Segre stratification matches the stratification given by the higher secant varieties: For a general $V \in \MS_X(2n;2e)$, let $E \in M(V)$ so that 
$E$ is general in $U_{X}(n, -e)^s$. Let $k$ be the smallest integer satisfying $ \Sec^{k} \pp E = \pp H^{1}(X, \Sym^{2}E)$. 
Then by Criterion \ref{geomlift}, there is some elementary transformation $F$ of $E^*$ with $\deg E/F = k$ lifting to $V$ as a Lagrangian subbundle. By deforming $V$ and $E$ inside  $\sy_e$, we may assume that $F$ is general in $U_X(n,e-k)^s$. Now consider the symplectic extension
\[ 0 \to F \to V \to F^* \to 0. \]
In the same way, $\delta(V) \in \pp H^{1}(X, \Sym^{2}F)$ belongs to $\Sec^{k}\pp F$, since the elementary transformation $E \to F^*$ lifts to $V$. Note that 
\[
\dim \Sec^k \pp F \leq k(n+1) -1 \ < \ (n+1)(k-e) + \frac{1}{2}n(n+1)(g-1) -1 = \dim \pp H^1 (X, \Sym^2 F),
\]
since we assumed $2e < n(g-1)$. Therefore, $\Sec^{k}\pp F$ is properly contained in  $\pp H^1 (X, \Sym^2 F)$, and certainly it is inside the closure of $\Sec^{k+1} \pp F \setminus \Sec^k \pp F$. Thus, again by Criterion \ref{geomlift}, the class $\delta(V)$ belongs to the closure of a family of bundles admitting liftings of elementary transformations of degree $(k-e)-(k+1) = -(e+1)$. In particular, $V$ belongs to the closure of $\MS_{X}(2n;2e+2)$ in $\MS_{X}(2n)$. By the irreducibility of $\MS_{X}(2n;2e)$, the same holds for arbitrary $V \in \MS_X(2n;2e)$. \qed 

\begin{remark} Suppose $t=2e < n(g-1)$. If $\gcd(n,e) = 1$, then $\tilde{U}_e = U_X(n,-e)^s$ and $\MS_X(2n;t)$ is birational to the fibration $\sy_e$ over $U_X(n,-e)^s$ whose fiber at $ E$ is $\pp H^1(X, \Sym^2 E)$. 
\end{remark}

Finally we give a geometric interpretation of the cardinality of $M(V)$, when it is finite.
Assume $t =2e = n (g-1) $. Let $V$ be general in $\MS_X(2n;t)$ and choose  $E \in M(V)$ with $\deg E = -e = -\frac{1}{2} n(g-1)$. Consider the subvariety $\pp E \subset \pp H^1(X, \Sym^2 E)$.
When $g \ge 4$, we have $\mu(E) < -1$, and by Lemma \ref{embedding}, the map $\phi_s : \pp E \to \pp H^1(X, \Sym^2 E)$ is an embedding. Furthermore, it was proven in \cite[Proposition 3.2]{CH2} that the secant variety $\pp E$ is not defective in $\pp H^1(X, \Sym^2 E)$ in the sense that
\[
\dim \Sec^{k} \pp E  \ = \ \min \{ (n+1)k -1, \ \dim \pp H^1 (X, \Sym^2 E) \}.
\]
In particular, for $k = n(g-1)$ we have 
\begin{equation}
\dim \Sec^{k} \pp E \ = \ (n+1)k -1 \  = \ \dim \pp H^1 (X, \Sym^2 E).
\label{nondefective} \end{equation}
By a \textsl{$k$-secant space of} $\pp E$, we mean a linear subspace in $\pp H^1(X, \Sym^2 E)$ spanned by some $k$ points of $\pp E$.
By (\ref{nondefective}), there are finite number  of $n(g-1)$-secant spaces of $\pp E$ which pass through a general point of $\pp H^1 (X, \Sym^2 E)$. The following result generalizes Lange-Narasimhan \cite[Proposition 2.4]{LN} for rank 2 bundles.
\begin{theorem} \ Assume $g \ge 4$ and $n(g-1)$ is even.
Let $V \in \MS_X(2n)$ be a general symplectic bundle with $t(V) = n(g-1)$, and let $E \in M(V)$.  
There is a one-to-one correspondence between $M(V) \setminus \{E\}$ and $n(g-1)$-secant spaces of $\pp E$ passing through the point $[E \subset V] \in \pp H^1(X, \Sym^2 E)$. \label{number}
\end{theorem}
\begin{proof}
By Criterion \ref{geomlift}, a general class $\delta(V)$ lies on a $n(g-1)$-secant space of $\pp E$ if and only if there is an associated elementary transformation $F$ of $E^*$ with $\deg(E^*/ F) = n(g-1)$, lifting to a Lagrangian subbundle of $V$. Since $t(V) = n(g-1)$, the Lagrangian subbundle $F$ is an element of $M(V)$. Conversely, by Theorem \ref{dimension} (3) every element of $M(V)$ other than $E$ appears in this way, since $V$ is general. 
 \end{proof}
\begin{remark}
It is an interesting problem to compute the cardinality of $M(V)$ explicitly. Let $N_{n,g}$ denote the cardinality of $M(V)$ for a general symplectic bundle $V$ of rank $2n$ over a curve of genus $g$ so that $n(g-1)$ is even. 
It is well known that $N_{1,g} = 2^g$. 
The same problem for maximal subbundles of vector bundles was solved by Holla \cite{Hol} (see also Lange--Newstead \cite{LN}). As far as we are aware, the number $N_{n,g}$ is not known in general.
\end{remark}

\section{Description of Segre strata for orthogonal bundles}

In this section, we investigate the geometry of the Segre stratification on the moduli space of orthogonal bundles $\MO_{X}(2n)$. In contrast to $\MS_{X}(2n)$, this has two connected components $\MO_{X}(2n)^{\pm}$, corresponding to bundles of trivial and nontrivial second Stiefel--Whitney class (see Serman \cite{Ser} for more details). We begin by determining the component to which each stratum $\MO_{X}(2n;t)$ belongs.

\subsection{Topological classification}

\begin{theorem} (1) Let $E_1$ and $E_2$ be isotropic rank $n$ subbundles of an orthogonal bundle $V$ of rank $2n$. Then $\deg E_1$ and $\deg E_2$ have the same parity. \\
(2) A semistable orthogonal bundle $V$ 
 belongs to $\MO_{X}(2n)^+$ (resp., $\MO_{X}(2n)^-$) if and only if its isotropic rank $n$ subbundles have even degree (resp., odd degree). \label{topology} \end{theorem}
\begin{proof} By Theorem \ref{stability}, the moduli maps $\alpha_{e} \colon \A_{e} \dashrightarrow \MO_{X}(2n)$ are defined on dense subsets for each $e>0$. Since each family $\A_e$ is connected, the image of each $\A_e$ is entirely contained in either $\MO_{X}(2n)^+$ or $\MO_{X}(2n)^-$. As was already noted in the proof of Theorem \ref{stability},  whenever $e_1$ and $e_2$ have the same parity, the images of $\A_{e_1}$ and $\A_{e_2}$  have nonempty intersection. Hence they 
lie on the same component of $\MO_X(2n)$. This is essentially due to the fact that any antisymmetric principal part has an even degree (see Lemma \ref{anti-irr}). 

  For each $k>0$, the image of $\A_{2k}$ is contained in $\MO_{X}(2n)^+$ since the trivial bundle has trivial second Stiefel--Whitney class. Thus the image of $\A_{2k-1}$ is contained in $\MO_{X}(2n)^-$. This shows the above statements (1) and (2) for the case when $V$ is a stable orthogonal bundle and the isotropic subbundles are stable bundles. By deforming an arbitrary bundle to a stable one, we  see that (1) and (2) hold in general.
\end{proof}

\subsection{An extra stratum}

The local deformations of an orthogonal bundle $V$ are given by $H^1(X, \wedge^2 V)$. Since $V$ is self-dual, $\wedge^2 V \subset \Endom(V)$. Hence both components $\MO_{X}(2n)^{\pm}$ have dimension
\[ h^{1}(X, \wedge^{2}V) = n(2n-1)(g-1), \]
since we may assume that $V$ is simple. On the other hand, 
\begin{eqnarray*}
\dim \A_e & = & \dim U_X(n,-e)  + \dim \pp H^1 (X, \wedge^2 E) \\
& = & (n-1)e + \frac{1}{2} n (3n-1) (g-1)
\end{eqnarray*}
for a general $E \in U_X(n,-e)$. Comparing dimensions, we see that in order for an $\A_e$ to cover either component of $\MO_{X}(2n)$, we need $e \geq \left\lceil \frac{1}{2}n(g-1) \right\rceil$. However, a connected set $\A_e$ can cover at most one component. Hence we are led to consider at least two top strata, corresponding to distinct values of $e \geq \left\lceil \frac{1}{2}n(g-1) \right\rceil$. This is one important difference between the stratifications in the symplectic and orthogonal cases.

\subsection{Geometry of the strata}
The arguments in this subsection are essentially the same as those appearing in \S4, except for the complication coming from the topological property described in Theorem \ref{topology}.

\begin{theorem} Consider a positive even integer $t = 2e \leq n(g-1)+3$.
\begin{enumerate}
\renewcommand{\labelenumi}{(\arabic{enumi})}
\item A general point of $\A_e$ corresponds to a bundle $V$ with $t(V) = t$; in particular, $\MO_{X}(2n;t)$ is nonempty and irreducible.
\item If $t < n(g-1)$ then $M(V)$ is a single point for general $V \in \MO_{X}(2n;t)$.
\end{enumerate} \label{dimensionorth} \end{theorem}

\begin{proof} Let $E$ be a general bundle in $U(n, -e)^s$ for $e \ge 1$. As before, we will bound the dimension of the locus in $H^{1}(X, \wedge^{2}E)$ of extensions $V$ which admit a rank $n$ isotropic subbundle $F$ of degree $\ge -e$ other than $E$. Again, either $F$ is an elementary transformation of $E^*$ lifting to $V$, or $F$ fits into a diagram of the form (\ref{nontransv}).

\textit{Step 1.} We first show that the latter situation does not arise for a general $V$ in $H^1(X, \wedge^2 E)$ for $t \le n(g-1) - 1$ (not $t \leq n(g-1) + 1$ as in the symplectic case; see Remark \ref{ft}). For each fixed $H \subset E$, by Criterion \ref{genlift} (2), the locus of extensions fitting into a diagram of the form (\ref{nontransv}) for some $F$ of degree $\ge -e$ is bounded by
\[ \dim \left( \Sec^{e-h}\Gr(2, E/H) \right) + 1 + \dim \Ker \left( q_{*} \, ^tq^{*} \right). \]
As in the symplectic case, we obtain $e-h \geq 0$. Since
\[ \dim \Gr(2,E/H) = 2(n-r-2) + 1 = 2(n-r) - 3, \]
the secant variety has dimension bounded by
\[ (e-h)(2(n-r) - 3) + (e-h) - 1 = 2(e-h)(n-r-1) - 1. \]
By Proposition \ref{qq}, we have
\begin{align*} \dim \Ker \left(q_{*} \, ^tq^{*} \right) &= h^{1}(X, \wedge^{2}E) - h^{1} \left(X, \wedge^{2}(E/H) \right) \\
 &\leq h^{1}(X, \wedge^{2}E) - (n-r-1)(e-h) - \frac{1}{2} (n-r)(n-r-1)(g-1). \end{align*}
As before, the subbundle $H$ of $E$ varies in a Quot scheme of dimension
\[ h^{0}(X, \Homom(H, E/H)) = nh - re - r(n-r)(g-1). \]
Thus the dimension of the locus of extensions in $H^1(X, \wedge^2 E)$ admitting a diagram of the form (\ref{nontransv}) for some $F$ of degree $\ge -e$ is bounded by
\begin{eqnarray*}
\dim \left( \Sec^{e-h} \Gr(2, E/H) \right) + 1 + \dim \Ker (q_* \, ^t q^*) + h^0(X, \Homom (H, E/H)) \\
\le \ (e-h)(n-r-1) - \frac{1}{2} (n-r)(n-r-1) (g-1) + nh-re + h^1(X, \wedge^2 E) \\
= \ e(n-r) - (e-h)(r+1) - \frac{1}{2} (n-r)(n-r-1) (g-1) + h^1(X, \wedge^2 E).
\end{eqnarray*}
Since $e-h \geq 0$ and we have assumed $e \le \frac{1}{2} (n (g-1)-1)$, this is smaller than $h^1(X, \wedge^2 E)$. Therefore, for $t = 2e \le  n (g-1) -1$, a general $V \in H^1(X, \wedge^2 E)$ does not admit a rank $n$ isotropic subbundle $F$ of degree $\ge -e$ which fits into a diagram of the form (\ref{nontransv}).

\textit{Step 2.} Next we consider the case $r = 0$, so $F$ is an elementary transformation of $E^*$ lifting to $V$ isotropically. Write $\deg F = -f \ge -e$. By Criterion \ref{cohomlift} (3) and Lemma \ref{anti-irr}, we have $e+f \equiv 0 \mod 2$. By Criterion \ref{geomlift} (2), the dimension of the locus of extensions in $H^1(X, \wedge^2 E)$ admitting such an $F$ is bounded by
\[
\dim \left( \Sec^{\frac{1}{2}(e+f)} \Gr(2, E) \right) + 1 \ \le \ \frac{1}{2}(e+f) (2(n-2)+2) = (e+f)(n-1).
\]
If $f < \frac{1}{2} n (g-1)$, this bound is smaller than
\[
 h^1(X, \wedge^2 E)  = (n-1)e + \frac{1}{2}n(n-1)(g-1)
\]
so a general extension in $H^{1}(X, \wedge^{2}E)$ does not admit any such lifting.

\textit{Step 3} of the proof of Theorem \ref{dimension} carries over, mutatis mutandis, to the orthogonal case, to prove statement (2), and (1) for $t \leq n(g-1)-1$. To prove (1) for $n(g-1) \leq t \leq n(g-1)+3$, we just repeat the arguments in Steps 1 and 2 under the assumption $-f \geq -e + 2$, to show that for general $E \in U_{X}(n, -e)$, a general $V \in H^{1}(X, \wedge^{2}E)$ does not admit a rank $n$ isotropic subbundle $F$ of degree $\geq -e + 2$.
\end{proof}

\begin{remark} Note that we do not prove an analogue of Theorem \ref{dimension} (3) for the orthogonal case. A dimension count analogous to that in the proof of Theorem \ref{dimension} (3) does not exclude the possibility that two maximal isotropic rank $n$ subbundles of a general bundle in $\MO_{X}(2n)$ intersect in a line bundle, or in a rank 2 subbundle if $g = 2$. It is unclear to us at this stage whether or not this in fact happens. \label{ft} \end{remark}

Next, we compute the dimension of $M(V)$ for a general $V \in \MO_X(2n)^{\pm}$. We will use the following analogue of Lemma \ref{Zariskitangent}, which is proven in the same way as in the symplectic case:
\begin{lemma} \label{tangentorth}
For an orthogonal bundle $V \in \MO_X(2n)$ and a point $[\eta: E \subset V]$ of $M(V)$, the Zariski tangent space of $M(V)$ at $\eta$ is identified with $H^0(X, \wedge^2 E^*)$. \qed
\end{lemma}

\noindent Combining with Lemma \ref{Hirlemma} (2), we obtain:
\begin{corollary}
Let $t = 2e = n(g-1) + \e$ with $\e \in \{0,1,2,3 \}$. Then for a general $V \in \MO_X(2n;t)$, we have
\[ \dim M(V) = \frac{\varepsilon(n-1)}{2}. \] \qed
\end{corollary}

\noindent {\it Proof of Theorem \ref{mainorth}.}\\
\\
(1) Consider the rational map $\alpha_e: \A_e \to \MO_X(2n)$.  For each $t = 2e = n(g-1) + \varepsilon$ where $\varepsilon \in \{0,1,2,3 \}$, a direct computation shows that
\[
\dim M(V) = \dim \A_e - \dim \MO_X(2n).
\]
Since $\dim M(V) \ge \dim \alpha_e^{-1} (V)$ for any $V \in \MO_X(2n;2e)$, this equality implies that $\alpha_e$ is dominant to a component of $\MO_X(2n)$. This shows that for general $V$ in some one of the components $ \MO_X(2n)^{\pm}$, we have
\[
n(g-1) \le t(V) \le n(g-1) +3.
\]
By semicontinuity, $t(V) \le n(g-1) +3$ for any orthogonal bundle $V$.

\noindent (2) The nonemptiness and irreducibility were proven in Theorem \ref{dimensionorth} (1). The remaining part can be proven in the same way as Theorem \ref{mainsymp} (2).

\noindent (3) This was proven in Theorem \ref{dimensionorth} (2). \qed

\subsection{Configuration of the dense strata} \label{config}

For each $2e = n(g-1 ) + \e$ with $0 \le \e \le 3$, the locations of the images of the $\alpha_{e}$ depend on the congruence class of $n(g-1)$ modulo 4: the trivial bundle of rank $2n$ is contained in $\MO_X(2n)^+$, hence for each $k$, we have $\MO_X(2n;4k) \subset \MO_X(2n)^+$ and $\MO_X(2n;4k+2) \subset \MO_X(2n)^-$.
We may summarize the situation for the dense strata as follows:
\[ n(g-1) \equiv 0 \mod 4: \quad \quad \begin{array}{c|c|c} t & \hbox{Component} & \dim M(V) \\ \hline
n(g-1) & \MO_{X}(2n)^{+} & 0 \\
n(g-1) + 2 & \MO_{X}(2n)^{-} & \hbox{$\frac{1}{2}(n-1)$} \end{array} \]

\[ n(g-1) \equiv 1 \mod 4: \quad \quad \begin{array}{c|c|c} t & \hbox{Component} & \dim M(V) \\ \hline
n(g-1) + 1 & \MO_{X}(2n)^{-} & \hbox{$\frac{1}{2}(n-1)$} \\
n(g-1) + 3 & \MO_{X}(2n)^{+} & \hbox{$n-1$} \end{array} \]

\[ n(g-1) \equiv 2 \mod 4: \quad \quad \begin{array}{c|c|c} t & \hbox{Component} & \dim M(V) \\ \hline
n(g-1) & \MO_{X}(2n)^{-} & \hbox{0} \\
n(g-1) + 2 & \MO_{X}(2n)^{+} & \hbox{$\frac{1}{2}(n-1)$} \end{array} \]

\[ n(g-1) \equiv 3 \mod 4: \quad \quad \begin{array}{c|c|c} t & \hbox{Component} & \dim M(V) \\ \hline
n(g-1) + 1 & \MO_{X}(2n)^{+} & \hbox{$\frac{1}{2}(n-1)$} \\
n(g-1) + 3 & \MO_{X}(2n)^{-} & \hbox{$n-1$} \end{array} \]

\begin{remark} In \cite[\S 4]{CH2}, symplectic bundles $W$ of rank four were studied which satisfy $s_{2}(W) < t(W)$; that is, whose maximal vector subbundles of half rank are all nonisotropic. Orthogonal bundles $V$ belonging to the top stratum $\MO_{X}(2n;n(g-1)+2)$ provide a different example of this phenomenon. Due to the Hirschowitz bound \cite{Hir}, all such $V$ have a vector subbundle of degree at least $- \left\lceil \frac{1}{2}n(g-1) \right\rceil$, but no isotropic rank $n$ subbundle of this degree or greater.

When $n(g-1)$ is even, such an orthogonal $V$ is nongeneral as a vector bundle (compare with Lange--Newstead \cite[\S 2]{LaNe}): Since any maximal rank $n$ subbundle $F \subset V$ is nonisotropic, we must have $h^{0}(X, \Sym^{2}F^{*}) > 0$. But in this case $\mu(\Sym^{2}F^{*}) \leq g-1$. By Lemma \ref{Hirlemma} (2), this is a condition of positive codimension on $F$, so none of the maximal vector subbundles of $V$ are general. \end{remark}

\noindent \footnotesize{Department of Mathematics, Konkuk University, 1 Hwayang-dong, Gwangjin-Gu, Seoul 143-701, Korea.\\
E-mail: \texttt{ischoe@konkuk.ac.kr}\\
\\
H\o gskolen i Oslo og Akershus, Postboks 4, St. Olavs plass, 0130 Oslo, Norway.\\
E-mail: \texttt{george.hitching@hioa.no}}

\end{document}